\documentclass[12 pt]{article}

\usepackage{verbatim}

\usepackage{amsmath}

\usepackage{latexsym}

\usepackage{amssymb}

\usepackage{amscd}

\newtheorem{Thm}{Theorem}[section]

\newtheorem{Prop}[Thm]{Proposition}

\newtheorem{Lem}[Thm]{Lemma}

\newtheorem{Cor}[Thm]{Corollary}

\newtheorem{Def}[Thm]{Definition}

\newtheorem{Ex}[Thm]{Example}

\begin{document}

\openup8pt

\title{Frequency Permutation Arrays}

\author{Sophie Huczynska and Gary L. Mullen}
\date{}
\maketitle

\begin{abstract}
Motivated by recent interest in permutation arrays, we introduce
and investigate the more general concept of frequency permutation
arrays (FPAs).  An FPA of length $n=m \lambda$ and distance $d$ is
a set $T$ of multipermutations on a multiset of $m$ symbols, each
repeated with frequency $\lambda$, such that the Hamming distance
between any distinct $x,y \in T$ is at least $d$.  Such arrays
have potential applications in powerline communication.  In this
paper, we establish basic properties of FPAs, and provide direct
constructions for FPAs using a range of combinatorial objects,
including polynomials over finite fields, combinatorial designs,
and codes.  We also provide recursive constructions, and give
bounds for the maximum size of such arrays.

\end{abstract}

\section{Introduction}
As indicated in \cite{ChuColDuk} and \cite{ColKloLin}, permutation
arrays arise in the study of permutation codes, which in turn have
a natural applicability to powerline communications. An electric
power line may be used to transmit information in addition to
electric power, by modulating its frequency to form a set of $n$
close frequencies. These small variations may be decoded as
symbols at the receiver. Steps must be taken to ensure that this
information transmission does not interfere with the line's
primary function of power transmission, and for this reason block
coding is used (codewords of fixed length). A code is a
\emph{constant composition code} if each codeword, of length $n$,
has precisely $r_i$ occurrences of the $i$-th symbol, where the
$r_i $ are positive integers satisfying $\sum_{i=1}^{m}r_i=n$.
(Here, the $i$-th symbol corresponds to the $i$-th frequency.)
Various tradeoffs must be made between the competing goals of
addressing noise problems and the requirement of a constant power
envelope.

One approach is to choose $r_1=r_2=\cdots =r_n =1$, in which case
each codeword is a permutation on $n$ symbols.  An $(n,d)$
permutation array, usually denoted by $PA(n,d)$, is a set of
permutations of $n$ symbols with the property that the Hamming
distance between any two distinct permutations in the set is at
least $d$. Permutation arrays are important not only in powerline
communications as described above; they have also been applied in
the design of block ciphers; see \cite{delColLin}.

In this paper, we introduce a generalization of permutation
arrays, which we call {\it frequency permutation arrays}.  These
arise from the constant composition codes when we take
$r_1=r_2=\cdots=r_m=\lambda$, for some $\lambda$ such that
$n=m\lambda$.  When $\lambda=1$, this reduces to the permutation
case studied in \cite{ChuColDuk} and \cite{ColKloLin}. We present
various results and constructions for frequency permutation
arrays, many of which have well-known permutation array results as
special cases.  There is a strong connection with recent work on
constant composition codes such as \cite{ChuColDuk2} and
\cite{DinYin}.

\section{Frequency permutation arrays}

We consider rearrangements of the $n$-element set $\{1,
\ldots,1,2, \ldots,2, \ldots,m, \ldots,m \}$, where $n=m\lambda$
($n,m,\lambda \in \mathbb{N}$) and each of the $m$ distinct
symbols occurs exactly $\lambda$ times.  When $\lambda=1$, the set
of all such permutations is the symmetric group $S_n$ of
permutations on $n$ symbols.  In the general case, these
rearrangements are \emph{multipermutations} on the multiset $\{
1,\ldots,1,2, \ldots,2, \ldots,m, \ldots,m \}$ (each symbol
occuring $\lambda$ times); we shall call them
\emph{$\lambda$-permutations}.

\begin{Def}
Two distinct $\lambda$-permutations $\sigma=s_1\ldots s_n$, $\tau=t_1\ldots t_n$ have \emph{distance} $d(\sigma,\tau)=d$ if they disagree in $d$ entries, i.e. if $|\{i:s_i \neq t_i \}|=d$.  \end{Def}

This is the Hamming distance familiar from coding theory.  In the case when $\lambda=1$, two permutations $\sigma,\tau \in S_n$ have distance $d$ if $\sigma \tau^{-1}$ has exactly $n-d$ fixed points.

\begin{Def}
A \emph{permutation array} of length $n$ and minimum distance $d$,
denoted by $PA(n,d)$, is a subset $T$ of $S_n$ such that the
distance between any two members of $T$ is at least $d$.  A
$PA(n,d)$ may be viewed as an $s \times n$ array whose rows are
the $s$ permutations of $T$ in image form; taken pairwise, any two
distinct rows differ in at least $d$ positions.  The maximum
possible size of a $PA(n,d)$ is denoted by $M(n,d)$.
\end{Def}

We define a frequency permutation array as follows.

\begin{Def}
A \emph{frequency permutation array} of length $n=m \lambda$,
frequency $\lambda$ and minimum distance $d$, denoted by
$FPA_\lambda(n,d)$, is a set $T$ of $\lambda$-permutations
(multipermutations of the multiset
$\{1,\ldots,1,2\ldots,2,\ldots,m,\ldots,m\}$, each symbol repeated
$\lambda$ times), with the property that the distance between any
two members of $T$ is at least $d$. Equivalently, an
$FPA_\lambda(n,d)$ is an $s\times n$ array whose $s$ rows consist
of $m$ distinct symbols, each repeated exactly $\lambda$ times,
such that the distance between any two rows is at least $d$.
\end{Def}

Thus an $FPA_1(n,d)$ is simply a $PA(n,d)$. We let $M_ \lambda
(n,d)$ denote the maximum possible number of rows that can exist
in any $FPA_\lambda (n,d)$; then $M_1(n,d)=M(n,d)$.

\begin{Ex}  An $FPA_{3}(6,4)$ of size $4$ is given by
$$
L_1= \begin{array}{cccccc}
1&1&1&0&0&0 \\
1&0&0&1&1&0 \\
0&1&0&1&0&1 \\
0&0&1&0&1&1
\end{array}
$$
\end{Ex}

We first establish basic properties of frequency permutation
arrays. Various basic results on permutation arrays (for example
from \cite{ChuColDuk} and \cite{ColKloLin}) appear as special
cases of these results.

\begin{Thm}\label{Basic}  Let $n=m \lambda$.  Then
\begin{itemize}
\item[(i)] $M_\lambda (n,2) = \frac {n!} {(\lambda !)^m}$;

\item[(ii)] $M_\lambda (n,n) = m$;

\item [(iii)] $M_\lambda (n,d) \geq M_\lambda (n-\lambda,d),
M_\lambda (n, d+1)$;

\item[(iv)] If $n_1=m \lambda_1$ and $n_2=m \lambda_2$, then
\[M_{\lambda_1+\lambda_2}(n_1+n_2,d_1+d_2) \geq
\mathrm{min}\{M_{\lambda_1}(n_1,d_1),M_{\lambda_2}(n_2,d_2)\}.\]
In particular, $M_{2\lambda}(2n,2d) \geq M_{\lambda}(n,d)$.

\item[(v)] $M_{\lambda}(n,d) \leq \frac{M(n,d)}{\lambda}$;
$M_{\lambda}(n,d) \leq \frac{n!}{\lambda(d-1)!}$.

Moreover, for any divisor $l$ of $\lambda$, $M_{\lambda}(n,d) \leq
\frac{l}{\lambda} M_l(n,d)$.
\end{itemize}
\end{Thm}

\noindent {\bf Proof} (i) Since two distinct multipermutations
must differ in at least two entries, $M_{\lambda}(n,2)$ is the
number of distinct $\lambda$-permutations.  There are
$\binom{n}{\lambda}.\binom{n- \lambda}{\lambda} \ldots
.\binom{n-(m-1)\lambda}{\lambda}$ choices for each
$\lambda$-permutation, i.e. $\binom{m
\lambda}{\lambda}.\binom{(m-1)\lambda}{\lambda} \ldots
.\binom{\lambda}{\lambda}=\frac{(m \lambda)!}{({\lambda}!)^m}$
such multipermutations in total.

(ii) Since there are at most $m$ choices for the symbol in the
first position of a $\lambda$-permutation in an
$FPA_{\lambda}(n,n)$, we have $M_{\lambda}(n,n) \leq m$.  Take $m$
blocks comprising $\lambda$ copies of each symbol: $\{0 \ldots
0\}$, $\{1, \ldots, 1\}$, $\ldots,$ $\{m-1, \ldots, m-1 \}$;
applying an $m$-cycle to these blocks yields $m$
$\lambda$-permutations, all of pairwise distance $n$.

(iii) Adding $\lambda$ copies of some new symbol to each row of an
$FPA_{\lambda}(n-\lambda,d)$ yields an $FPA_{\lambda}(n,d)$; the
second observation is immediate from the definition.

(iv) Juxtaposing an $FPA_{\lambda_1}(n_1,d_1)$ and an
$FPA_{\lambda_2}(n_2,d_2)$ yields an $FPA_{\lambda}(n,d)$ with
$\lambda=\lambda_1+\lambda_2$, $n=n_1+n_2$ and $d=d_1+d_2$.

(v) This is proved in Theorem \ref{UpperBound}; from
\cite{ChuColDuk}, the size of a $PA(n,d)$ is bounded above by
$\frac{n!}{(d-1)!}$. $\hfill \square$

For any $\lambda$-permutation $\sigma$, the sphere with centre
$\sigma$ and radius $r$ is defined to be the set of all
$\lambda$-permutations with distance at most $r$ from $\sigma$.
We denote its volume by $V_{\lambda}(n,r)$.

\begin{Lem} Let $n=m \lambda$.  Then
\[ V_{\lambda}(n,r) = 1+\sum_{k=1}^{r} \sum_{P(k)} \frac{m!}{{r_1}! \ldots {r_s}! (m-t)!}
\binom{\lambda}{k_1}\binom{\lambda}{k_2}\ldots\binom{\lambda}{k_t}(-1)^k
\int_{0}^{\infty} e^{-x} \{ \prod_{j=1}^{t}L_{k_j }(x)\} dx,\]
where $L_k(x)$ is the $k$th Laguerre polynomial.  Here the inner
sum runs over $P(k)=\{(k_1,\ldots,k_t;r_1,\ldots,r_s)\}$, the set
of all partitions $k_1+ \cdots + k_t$ of $k \in \mathbb{N}$ into
positive integers $1 \leq k_i \leq \lambda$, where the set
$\{k_1,\ldots,k_t\}$ consists of $r_j$ occurrences of value
$k_{i_j}$ ($j=1,\ldots,s$), with $1 \leq k_{i_j} \leq \lambda$, $1
\leq r_j \leq t$ and $r_1+ \cdots + r_s=t$.
\end{Lem}

\noindent {\bf Proof}  Let $\sigma$ be any $\lambda$-permutation
of length $n$.  The set of $\lambda$-permutations at distance $k$
from $\sigma$ is obtained by taking each $k$-entry subset of
$\sigma$, and deranging its entries.  By a result obtained in
\cite{EveGil}, and reproved in \cite{Car}, the number of
derangements of a sequence composed of $n_1$ objects of type $1$,
$n_2$ objects of type $2$,$\ldots$, $n_t$ objects of type $t$
(i.e. permutations in which no object occupies a site originally
occupied by an object of the same type) is given by
\[ D(n_1,\ldots,n_t)=(-1)^N \int_{0}^{\infty} e^{-x} \{
\prod_{j=1}^{t} L_{n_j}(x)\} dx,\] where $n_1+\cdots+n_t=N$. The
result follows upon applying this theorem to each $k$-element
subset of $\sigma$. For any $\lambda$-permutation $\sigma$, and
any partition $k_1+ \cdots + k_t$ of $k \in \mathbb{N}$ into
positive integers $1 \leq k_i \leq \lambda$ ($1 \leq t \leq m$),
we count the number of $k$-subsets comprising $k_1$ occurrences of
symbol $s_1$, $k_2$ occurrences of symbol $s_2$, $\ldots$, $k_t$
occurrences of symbol $s_t$. Suppose the set $\{k_1,\ldots,k_t\}$
consists of $r_1$ occurrences of value $k_{i_1}$, $\ldots$, $r_s$
occurrences of value $k_{i_s}$, where $r_1+ \cdots + r_s=t$. There
are $\binom{m}{r_1} \binom{m-r_1}{r_2} \cdots
\binom{m-{\sum_{i=1}^{s-1}r_i}}{r_s}=\frac{m!}{r_1 ! \ldots r_s !
(m-t)!}$ choices for symbols $s_1,\ldots,s_t$.  For each choice,
there are
$\binom{\lambda}{k_1}\binom{\lambda}{k_2}\ldots\binom{\lambda}{k_t}$
subsets of $\sigma$ in which elements occur with appropriate
frequency. $\hfill \square$

A covering argument yields the following lower bound for
$M_{\lambda}(n,d)$, an analogue of the Gilbert-Varshamov bound in
coding theory, while a sphere-packing argument yields an upper
bound, analogous to the Hamming bound for coding.
\begin{Thm}
\[ \frac{n!}{(\lambda !)^m V_{\lambda}(n,d-1)} \leq
M_{\lambda}(n,d) \leq \frac{n!}{(\lambda !)^m V_{\lambda}(n,
\lfloor \frac{d-1}{2} \rfloor) }.\]
\end{Thm}

We remark in passing that a useful upper bound for the maximum
size of general constant-composition codes (CCCs) has recently
been presented in \cite{LuoFuHanChe} and has been further
developed in \cite{DinYin}.  However, for a CCC in which all
symbols of a codeword occur with equal frequency $\lambda$ (the
situation corresponding to FPAs), this upper bound essentially
reduces to the Plotkin bound, $M_{\lambda}(n,d) \leq
\frac{d}{d-n+\lambda}$, which is valid only when $d>n-\lambda$.
Since the direct constructions presented in this paper have
minimum distance less than or equal to $n-\lambda$, the bound is
of limited applicability in our setting.

\section{Direct constructions}

It is known that permutation arrays may be constructed using latin
squares (see \cite{ChuColDuk} and \cite{HucMul}).  Frequency
permutation arrays are related to frequency squares as permutation
arrays are to latin squares, and this connection may be exploited
to obtain a construction for FPAs.

Recall that a \emph{latin square of order $n$} is an $n \times n$
array in which $n$ distinct symbols are arranged so that each
symbol occurs once in each row and column.  Two latin squares
$L_1$ and $L_2$ of the same  order $n$ are said to be
\emph{orthogonal} if, when superimposed, each of the possible
$n^2$ ordered pairs occurs exactly once.  A set $\{ L_1, L_2,
\ldots, L_t \}$ of $t\geq 2$ latin squares is said to be
\emph{mutually orthogonal} (a set of MOLS) if the squares in the
set are pairwise orthogonal. Latin squares have been generalized
to allow repetitions of elements in each row and column.

\begin{Def}Let $n=m \lambda$. An $F(n;\lambda)$ {\it frequency square}
is an $n \times n$ array in which each of $m$ distinct symbols
occurs exactly $\lambda $ times in each row and column. Moreover
two such squares are {\it orthogonal} if when superimposed, each
of the $m^2$ possible ordered pairs occurs $\lambda ^2$ times.
\end{Def}

The following result in fact contains Proposition 1.2 of
\cite{ChuColDuk} as a special case.

\begin{Thm} \label{MOFS}
If there are $E$ mutually orthogonal frequency squares of type
$F(n;\lambda)$ where $n=m \lambda$, then $M_\lambda (n \lambda, n
\lambda -\lambda ^2) \geq mE$.  In particular, if $q$ is a prime
power and $i\geq 1$ is a positive integer, then \[ M_{q^{i-1}}
(q^{2i-1}, q^{2i-1}-q^{2i-2}) \geq q(q^i -1)^2/(q-1).\] Further if
$i=1$, $M_1(q,q-1)=q(q-1)$.
\end{Thm}

\noindent {\bf Proof} Label the rows and columns of each $n \times
n$ frequency square by the elements $0,1, \dots, n-1$. Then from
each of the frequency squares, form a set of $n \lambda$-tuples as
follows. For each symbol $i=0,1, \dots, m-1$, form an $n
\lambda$-tuple by listing the cell locations $(k,l)$ where $i$
occurs in the given square, proceeding row-by-row as $k$ runs from
$0$ to $n-1$.  Viewed as $m$ blocks, each of size $n\lambda$, of
an affine resolvable design, these form a parallel class of size
$m$. In total from the $E$ squares, $Em$ such $n \lambda$-tuples
are obtained, corresponding to $E$ parallel classes. The entries
of each $n \lambda$-tuple are ordered pairs; form new $n
\lambda$-tuples by disregarding the first coordinate of each
ordered pair. The resulting $n \lambda$-tuples form the rows of an
$FPA_\lambda(n \lambda, n\lambda -\lambda ^2)$. For, since each
symbol occurs $\lambda$ times in each column of a frequency
square, each row of the array comprises $\lambda$ copies of each
of the $n$ column-headings. Any two rows of the FPA arising from
the same parallel class will have distance $n \lambda$. Any two
rows derived from different classes will, due to the orthogonality
of the corresponding frequency squares, agree in at most $\lambda
^2$ positions, since agreement in $p$ positions implies that some
ordered pair occurs $p$ times when the MOFS are juxtaposed.  Hence
the array has minimum distance $n \lambda -\lambda ^2$. $\hfill
\square$

For $n=m \lambda$, it is known that the maximum number of mutually
orthogonal frequency squares (MOFS) of the form $F(n;\lambda)$ is
bounded above by $(n-1)^2/(m-1)$. Further, if $q$ is any prime
power and $i\geq 1 $ is a positive integer, then using linear
polynomials in $2i$ variables over the finite field $F_q$, a
complete set of $F(q^i;q^{i-1})$ MOFS can be constructed.
Specifically, take the polynomials $a_1x_1 + \cdots +
a_{2i}x_{2i}$ where neither $(a_1, \dots, a_i)$ nor $(a_{i+1},
\dots, a_{2i})$ is the zero vector $(0, \dots, 0)$ and no two of
the vectors are nonzero $F_q$ multiples of each other, i.e.
$({a_1}^\prime, \dots, {a_{2i}}^\prime) \not = e(a_1, \dots,
a_{2i})$ for any nonzero $e\in F_q$.  Further details may be found
in Chapter 4 of \cite{LayMul}.

We remark in passing that, while the array obtained from Theorem
\ref{MOFS} is optimal in size when $i=1$, it is not necessarily
optimal for $i>1$. This is in some sense expected because, in
using these complete sets of mutually orthogonal frequency squares
to construct error-correcting codes, the resulting codes are
maximal distance separable only in the case when $i=1$; see
\cite{DenMulSuc}. For example in the case $q=i=2$, Theorem
\ref{SepPA} yields an $FPA_2(8,4) $ with more than 18 rows (see
Example \ref{FPA284}).

Another way to build frequency permutation arrays utilises finite
fields, and may be considered as extending the approach of Theorem
2.4 of \cite{ChuColDuk}.

\begin{Thm}\label{LinPoly}
Let $L(x)=\sum_{s=0}^{i-1} \alpha_s x^{q^s} \in F_{q^i}[x]$.
Denote by $q^l$ the degree of $L(x)$, and by $r$ the rank of the
matrix
\[ A(L)= \left(
\begin{array}{ccccc}
\alpha_0 & \alpha_{i-1}^q & \alpha_{i-2}^{q^2} & \ldots & \alpha_1^{q^{i-1}}\\
\alpha_1 & \alpha_0^q & \alpha_{i-1}^{q^2} & \ldots & \alpha_2^{q^{i-1}}\\
\alpha_2 & \alpha_1^q & \alpha_0^{q^2} & \ldots & \alpha_3^{q^{i-1}}\\
\vdots & \vdots & \vdots & \vdots & \vdots \\
\alpha_{i-1} & \alpha_{i-2}^q & \alpha_{i-3}^{q^2} & \ldots & \alpha_0^{q^{i-1}}\\
\end{array} \right),
\]
so that $1 \leq r \leq i$. Let $0< d < q^{i-l}$. Then
\[M_{q^{i-r}}(q^i,q^i-dq^l) \geq \sum_{j=1}^{d} \frac
{N_j(q^i)} {q^{i-r}},\] where $N_j(q^i)$ denotes the number of
permutation polynomials over $F_{q^i}$ of degree $j$.
\end{Thm}
\noindent {\bf Proof}  It is a well-known result (see p 361 of
\cite{LidNie}) that the linearized polynomial $L(x)$ is a
permutation polynomial of $F_{q^i}$ if and only if the determinant
of the matrix $A(L)$ is non-zero.  More generally, the value set
of $L$ has cardinality $q^r$, where $r$ is the rank of $A(L)$. So
the linear transformation on $F_{q^i}$ defined by the polynomial
$L(x)$ has image of cardinality $q^r$ and kernel of cardinality
$q^{i-r}$.  Note that $q^{i-r} \leq q^l$.

Form an array as follows: for each permutation polynomial $f(x)$
over $F_{q^i}$, form a row by taking the images of the function
$L(f(x))$ as $x$ runs through the elements of the field $F_{q^i}$.
Each row is a $\lambda$-permutation of length $q^i$ on $m=q^r$
symbols, each occurring with frequency $\lambda=q^{i-r}$. If
$f(x)$ and $g(x)$ are permutation polynomials over $F_{q^i}$ of
degrees at most $d$, then the polynomial $L(f(x))-L(g(x))$ has
degree at most $dq^l$.  Hence (unless it is the zero polynomial)
it has at most $dq^l$ roots in $F_{q^i}$, and so appropriately
chosen $f(x)$ and $g(x)$ yield distinct rows of distance at least
$q^i-dq^l$. We must now ensure that $L(f-g)$ is not the zero
polynomial. This happens if and only if the value set of the
polynomial $f-g$ lies wholly within the kernel of $L$, which has
cardinality $q^{i-r}$. Suppose first that $f-g$ is non-constant.
Now, $f-g$ has degree at most $d<q^{i-l}$ and, since a polynomial
of degree $d$ cannot have more than $d$ roots in a field, its
value set has cardinality at least $\lfloor \frac{q^i-1}{d}
\rfloor + 1
> q^l \geq q^{i-r}$.  So the value set of $f-g$ cannot be
contained entirely within the set of $q^{i-r}$ values mapped by
$L$ to zero, and hence $L(f-g)$ is not the zero polynomial. For
the constant case note that, for any permutation polynomial
$f(x)$, all $f(x)+c$ with $c \in F_{q^i}$ are also permutation
polynomials.  For $f(x)+c$ to yield distinct rows, $c$ must run
through precisely one representative for each coset of the kernel
of $L$; there are $q^r$ of these. Taking $\frac{q^r}{q^i}$ of the
total number of permutation polynomials yields the desired number
of rows. $\hfill \square$

Observe that, in Theorem \ref{LinPoly}, if we take $L$ to be the
permutation polynomial $x^{q^{i-1}}$, we have maximal rank $r=i$
and degree $q^l=q^{i-1}$, so we obtain a $PA(q^i,q^i-dq^{i-1})$ of
size $\sum_{j=1}^{d} N_j(q^i)$.

To build an FPA with desired parameters, appropriate linearized
polynomials may be chosen, whose properties are known in advance.
The following corollaries illustrate two such constructions.
Recall that the trace function $TR:F_{q^i} \rightarrow F_q$ is
defined for $\alpha \in F_{q^i}$ by $TR(\alpha) = \alpha +\alpha^q
+\alpha^{q^2} + \dots + \alpha^{q^{i-1}}$.  More generally,
letting $i=gh$ and setting $E=F_{q^i}$ and $F=F_{q^h}$, the trace
function $TR_{E/F}:E \rightarrow F$ is defined for $\alpha \in E$
by $$TR_{E/F}(\alpha)= \alpha + \alpha^{q^h}+ \alpha^{q^{2h}}+...+
\alpha^{q^{(g-1)h}}.$$

\begin{Cor}\label{Trace}
Let $q$ be a prime power and let $i$ and $h$ be positive integers
such that $h$ divides $i$. Let $0 <d < q^h$.  Then \[
M_{q^{i-h}}(q^i, q^i - dq^{i-h}) \geq \sum_{j=1}^{d} \frac
{N_j(q^i)} {q^{i-h}}.\]
\end{Cor}
\noindent {\bf Proof}  Let $i=gh$, let $E=F_{q^i}$ and
$F=F_{q^h}$. Take $L$ to be the generalized trace function
$TR_{E/F}$ defined above; its kernel has cardinality $q^{i-h}$ and
its value set has cardinality $q^h$.  For any two permutation
polynomials $f,g$ over $F_{q^i}$ of degree at most $d<q^h$, the
value set of (non-constant) $f-g$ has cardinality at least
$\lfloor \frac{q^i-1}{d} \rfloor +1
> q^{i-h}$ and so $TR(f-g)$ is not the zero polynomial.  As in the
proof of Theorem \ref{LinPoly} dividing by $q^{i-h}$ deals with
the case when $f-g$ is constant. $\hfill \square$

In the next section we consider how a PA may be converted into an
FPA by appropriate substitutions on its symbols. If $q$ is a prime
power, a natural choice might be to apply the trace function to
the rows of a $PA(q^i, d)$. However if, for example, there are two
rows in the $PA(q^i,d)$ which differ by a constant $a \in F_{q^i}$
with $TR(a)=0$, then the resulting two rows in the $FPA$ will be
identical and so the rows will have distance 0. Thus applying the
trace function to the elements of an arbitrary PA does not appear
to be a good method to apply in a general setting.

\begin{Cor}\label{LinPoly2}
Let $q$ be a prime power and let $i$ and $n$ be positive integers
such that $n(<i)$ divides $i$. Let $0< d < q^{i-n}$. Then
\[M_{q^n}(q^i,q^i-dq^n) \geq \sum_{j=1}^{d} \frac
{N_j(q^i)} {q^n},\] where $N_j(q^i)$ denotes the number of
permutation polynomials over $F_{q^i}$ of degree $j$.
\end{Cor}
\noindent {\bf Proof}  Let $L(x)=x^{q^n}-x$ in Theorem
\ref{LinPoly}; its roots are precisely the elements of $F_{q^n}$.
The polynomial $L$ defines a linear transformation on $F_{q^i}$
whose kernel is the subfield $F_{q^n}$ and whose value set has
cardinality $q^{i-n}$. For permutation polynomials $f,g$ of degree
at most $d<q^{i-n}$, non-constant $f-g$ has value set of
cardinality at least $\lfloor \frac{q^i-1}{d} \rfloor + 1 > q^n$
and so identical rows can arise only in the case when
$g(x)=f(x)+c$ with $c \in F_{q^n}$.$\hfill \square$

We refer to \cite{Das} for a method for computing the value of
$N_j(q^i)$ for any prime power $q$ and  positive integers $i$ and
$j$. Note, however, that the result of \cite{Das} requires
considerable computation to compute, and that the corresponding
permutation polynomials which arise from solutions to the system
of equations in \cite{Das} must be constructed before the FPA can
be built.

We next indicate how affine resolvable designs can be used to
construct FPAs. A \emph{balanced incomplete block design} consists
of a finite set $V$ of $v$ points, and a collection $B$ of equally
sized subsets of $V$ called blocks, each of size $k$, such that
every pair of distinct points of $V$ occurs in exactly $\lambda$
blocks.  A \emph{resolvable} design has the additional property
that the collection $B$ of blocks can be partitioned into
\emph{parallel classes} (or resolution classes), such that every
point of $V$ occurs exactly once in each parallel class.  An
\emph{affine resolvable design}(ARD) is a resolvable design with
the further property that any two non-parallel blocks intersect in
precisely $\mu$ points, where $\mu=\frac{k^2}{v} \in \mathbb{N}$.
When $\mu=1$, the ARD is an affine plane of order $k^2$.  Given an
ARD, by labelling the blocks in each class and listing the blocks
in which each point lies, an FPA may be constructed.
\begin{Thm}\label{ARD}
\begin{itemize}
\item[(i)] Given an affine resolvable $(v,k,\lambda)$ design with
$r$ parallel classes, an $FPA_k(v,v-k)$ may be constructed of size
$r$. \item[(ii)] If there exist $m$ MOLS of order $n$, then an
$FPA_n(n^2,n^2-n)$ may be constructed of size $m+2$.  In
particular, if $q$ is a prime power, an $FPA_q(q^2,q^2-q)$ may be
constructed of size $q+1$.
\end{itemize}
\end{Thm}

More details of this approach, including a proof of Theorem
\ref{ARD}(i), may be found in \cite{ChuColDuk2}. In the MOLS case,
a standard construction may be used to build $m+2$ parallel
classes of an affine plane from the $m$ MOLS, and these classes
then used in part (i) to form an $FPA_n(n^2,n^2-n)$. Equivalently,
this FPA may be constructed directly by writing the rows of each
of the $m+2$ latin squares side-by-side to form new rows of length
$n^2$; it is clear from latin square properties that the resulting
array is an $FPA_n(n^2,n^2-n)$.

\begin{Ex}\label{OAEx}
Using the following $2$ MOLS of order $3$
$$
L_1= \begin{array}{ccc}
0&1&2 \\
1&2&0 \\
2&0&1 \\
\end{array}, L_2= \begin{array}{ccc}
0&1&2 \\
2&0&1 \\
1&2&0 \\
\end{array}
$$
yields the following $FPA_3(9,6)$:
\[
\begin{array}{ccccccccc}
0&0&0&1&1&1&2&2&2\\
0&1&2&0&1&2&0&1&2\\
0&1&2&1&2&0&2&0&1\\
0&1&2&2&0&1&1&2&0.
\end{array}
\]
\end{Ex}

\begin{Def}
An \emph{orthogonal array} of size $v$, with $r$ constraints, $s$
levels and strength $t$, denoted $OA[v,r,s,t]$, is an $r \times v$
array with entries from a set of $s \geq 2$ symbols, having the
property that in every $t \times v$ submatrix, every $t \times 1$
column vector appears the same number $\frac{v}{s^t}$ of times.
\end{Def}

The frequency array constructed in Theorem \ref{ARD} is in fact an
orthogonal array of strength $2$.  This gives rise to the
following observation.

\begin{Prop}\label{OA}
Every orthogonal array $OA[v,r,s,2]$ of strength $2$ is an
$FPA_{\frac{v}{s}}(v,v-\frac{v}{s})$ of size $r$.
\end{Prop}
\noindent {\bf Proof} In any row, each of the $s$ symbols occurs
with frequency $\frac{v}{s}$. For any pair of rows, each of the
$s^2$ pairs $(i,j)$ of elements occurs $\frac{v}{s^2}$ times. In
particular, each of the $s$ pairs $(i,i)$ occurs $\frac{v}{s^2}$
times, and hence two rows agree pairwise in precisely
$\frac{v}{s}$ positions. $\square$

Note that the FPAs obtained in this way are \emph{equidistant}, in
the sense that any two rows have distance precisely
$v-\frac{v}{s}$.  For constructions of orthogonal arrays, see for
example \cite{ColDin}; their connection with affine resolvable
designs is explored in \cite{BaiMonMor}.

We end the section with a construction of frequency arrays from
MDS codes.  Recall that a $q$-ary $(n,k)$ code is said to be
maximal distance separable (MDS) if it satisfies the Singleton
bound with equality, i.e. if $d=n-k+1$.

\begin{Thm}
Given an $[n,k,d]$ MDS linear code $C$ over $F_q$, the array
formed by taking the codewords of $C$ as columns is an
$FPA_{q^{k-1}}(q^k, q^{k-1}(q-1))$.
\end{Thm}

\noindent {\bf Proof}  Let $C$ be an $[n,k,d]$ MDS linear code
over $F_q$.  Let $G$ be a $k \times n$ generator matrix for the
code $C$, and write $G=[C_1 C_2 \ldots C_n]$, where the $C_i$ are
the columns of $G$. Form an $n \times q^k$ array $A$ by taking the
codewords of $C$ as the columns of $A$. These are given by $G^T
x^T=C^T$ as $x$ runs through $F_q^k$; the rows
${C_1}^T,\ldots,{C_n}^T$ of $G^T$ can be viewed as generating the
rows of $A$.

Each element of $F_q$ occurs in each row of $A$ with frequency
$q^{k-1}$, i.e. occurs $q^{k-1}$ times as the $i$th coordinate of
the codewords of $C$. Let $g_1, ..., g_k$ be the elements in the
$i$-th column $C_i$ of $G$, and consider the equation $a_1g_1 +
... + a_kg_k = b$, where $b, a_1, \ldots,a_k \in F_q$. Since $C_i$
has at least one nonzero value, say in the $j$-th row, we can
isolate the term $a_jg_j =b- \sum_{l \neq j}a_l g_l$. Then we can
arbitrarily assign $q$ values to each of $k-1$ remaining $a$'s,
and uniquely solve the equation for $a_j$ since $g_j \not = 0$.
Thus there are $q^{k-1}$ solutions for each value of $b$ in the
$i$-th coordinate.

Consider the distance between the two rows of the FPA
corresponding to ${C_i}^T$ and ${C_j}^T$. We have the system of
equations ${C_i}^T \cdot(x_1,\ldots,x_k)= \alpha$ and ${C_j}^T
\cdot(x_1,\ldots,x_k) = \beta$ ($\alpha, \beta \in F_q$). For an
MDS code, any $k$ columns (in particular, any two columns) of the
generator matrix are linearly independent. Since ${C_i}^T$ and
${C_j}^T$ are linearly independent, this system of two linear
equations in $k$ variables will have rank $2$, and thus $q^{k-2}$
solutions. This means that every ordered pair $(\alpha,\beta)$
occurs $q^{k-2}$ times. Thus, in particular, the $q$ ordered pairs
$(\alpha,\alpha),\alpha \in F_q$ are obtained $q^{k-2}$ times, so
$A$ is an FPA with distance $q^k- qq^{k-2} =q^k-q^{k-1} =
q^{k-1}(q-1)$. $\hfill \square$

Observe that this provides a direct construction for a class of
FPAs described in Proposition \ref{OA}, with $v=q^k$ and $s=q$.
Example \ref{OAEx} may alternatively be obtained by the MDS
construction using the generator matrix
\[ G= \left(
\begin{array}{cccc}
1&0&1&2\\
0&1&1&1\\
\end{array} \right).
\]

\section{Constructing new FPAs from old}

In this section, we explore how one or more FPAs may be used as
ingredients in the construction of new FPAs.

\begin{Thm}\label{UpperBound}
\begin{itemize}
\item[(i)] Given an $FPA_{\lambda}(n,d)$ of size $N$, a $PA(n,d)$
may be constructed of size $\lambda N$.  In particular,
$M_{\lambda}(n,d) \leq \frac{M(n,d)}{\lambda}$.

\item[(ii)]  Let $l$ divide $\lambda$.  Given an
$FPA_{\lambda}(n,d)$ of size $N$, an $FPA_l(n,d)$ may be
constructed, of size $\frac{\lambda}{l}N$.  In particular,
$M_{\lambda}(n,d) \leq \frac{l}{\lambda}M_l(n,d)$.
\end{itemize}
\end{Thm}
\noindent {\bf Proof} (i) Denote the $FPA_{\lambda}(n,d)$ by $A$;
let the symbol set of $A$ be $\{0,1,\ldots,m-1\}$. Using
appropriate substitutions, $A$ can be converted to a $PA(n,d)$,
$A^{\prime}$, of size $N$.  For a row $R$ of $A$, moving from left
to right, replace the $\lambda$ occurrences of a given symbol $s$
by the sequence $s \lambda +1, s \lambda +2, \ldots, (s+1)\lambda$
($0 \leq s \leq m-1$).  The new row $R^{\prime}$ is a permutation
of $1,2,\ldots,n$.  Since agreement between any two rows of
$A^{\prime}$ can occur only at positions of agreement between the
corresponding rows of $A$, the PA $A^{\prime}$ has minimal
distance $d$.  Now perform a cyclic shift on the entries of each
substitution set $\{s \lambda +1, s \lambda +2, \ldots,
(s+1)\lambda \}$ ($0 \leq s \leq m-1$).  This process can be
repeated $\lambda$ times, to obtain $\lambda$ different
substitutions for $R$; all have pairwise distance $n$. Apply this
process to each row of $A$; the distance between new rows
corresponding to different rows of $A$ is at least $d$. Hence we
have a $PA(n,d)$ of size $\lambda N$.

(ii) The proof is analogous to that of part (i). In this case, the
substitution set for a given symbol $s$ of the
$FPA_{\lambda}(n,d)$ comprises $l$ copies each of $a$ symbols. The
generalization of the $\lambda$ cyclic shifts applied to the
substitution sets, is the set of $\frac{\lambda}{l}$ permutations
comprising an $FPA_l(\lambda,\lambda)$, described in part (ii) of
Theorem \ref{Basic}. $\hfill \square$

\begin{Ex}
An $FPA_6(12,6)$ is constructed in Example \ref{HadMat}. By
Theorem \ref{UpperBound}, this FPA can be converted first to an
$FPA_3(12,6)$ and then to a $PA(12,6)$.  We illustrate the use of
the substitutions (without the cyclic shifts) on four sample rows.
\end{Ex}
The first four rows of the $FPA_6(12,6)$ are
$$\begin{array}{cccccccccccc}
1&0&1&0&1&1&1&0&0&0&1&0\\
1&0&0&1&0&1&1&1&0&0&0&1\\
1&1&0&0&1&0&1&1&1&0&0&0\\
1&0&1&0&0&1&0&1&1&1&0&0
\end{array}
$$
After substitutions, four rows of the $FPA_3(12,6)$ are
$$\begin{array}{cccccccccccc}
3&1&3&1&3&4&4&1&2&2&4&2\\
3&1&1&3&1&3&4&4&2&2&2&4\\
3&3&1&1&3&1&4&4&4&2&2&2\\
3&1&3&1&1&3&2&4&4&4&2&2
\end{array}
$$
After substitutions, four rows of the $FPA_6(12,6)$ are
$$\begin{array}{cccccccccccc}
7&1&8&2&9&10&11&3&4&5&12&6\\
7&1&2&8&3&9&10&11&4&5&6&12\\
7&8&1&2&9&3&10&11&12&4&5&6\\
7&1&8&2&3&9&4&10&11&12&5&6
\end{array}
$$

Converting a PA to an FPA by substitution is less straightforward
in general.  The next result applies, for example, to an FPA
arising from an orthogonal array.

\begin{Prop}
Let $n=m \lambda$.  Let $A$ be an $FPA_{\lambda}(n,d)$ such that,
between any two rows, each of the $m^2$ pairs $(i,j)$ occurs
precisely $t$ times.  Then $A$ may be converted, by reduction mod
$r$ (where $r|m$) to an $FPA_{\frac{n}{r}}(n,n-\frac{tm^2}{r})$.
\end{Prop}
\noindent {\bf Proof}  Reduce the entries of $A$ mod $r$.  Each
row of the new array is a $\lambda$-permutation on $r$ symbols
with frequency $\frac{n}{r}$.  For any two rows in the new FPA,
the pair of entries $(i\, \mathrm{mod}\,r, j \,\mathrm{mod}\,r)$
agree, for each of the $\frac{m}{r}$ values of $j$ in the
congruence class of $i$. This yields $\frac{tm}{r}$ pairs for a
given value of $i$, yielding $\frac{tm^2}{r}$ such pairs in total,
i.e. a minimal distance of $n-\frac{tm^2}{r}$. $\hfill \square$

The substitution technique may also be used on permutation arrays
which have been constructed from latin squares.  For example,
given a PA obtained from Theorem \ref{MOFS} using a set of $q^i -
1$ MOLS of order $q^i$, applying the $(q^i -1 )/(q-1)$
substitutions from Theorem 9.20 of \cite{LayMul} to its entries,
yields an FPA as described in the second part of Theorem
\ref{MOFS}.  This approach allows the FPA to be built without
constructing the corresponding sets of MOFS.

In the proof of Theorem \ref{UpperBound}, a pairwise distance of
$n$ is imposed on the set of new rows derived from any given
original row.  Relaxing this condition to minimum distance $d$,
the $\lambda$-cycle (or its frequency analogue) may be replaced by
an appropriate (frequency) permutation array. This observation
underlies the next result.

\begin{Thm}
Let $n=m \lambda$, and let $F_1,\ldots,F_b$ be $b$
$FPA_{\lambda}(n,d)$'s (not necessarily different).  Let $C$ be an
$FPA_{n}(bn,c)$ of size $N$, where $c \geq bd$.  Then an
$FPA_{\lambda}(bn,bd)$ may be constructed, of size $N
\mathrm{min}_{1 \leq i \leq b}|F_i|$.
\end{Thm}

\noindent {\bf Proof}  Relabelling if necessary, construct
$F_1,\ldots,F_b$ on disjoint symbol sets, so there are $bn$
symbols in total. For each row in $C$, use the entries of the row
as column headings, and place the $n$ columns of each $F_i$ under
the $n$ occurrences of the symbol $i$. The resulting array is an
$FPA_{\lambda}(bn,bd)$, of size $\mathrm{min}_{1 \leq i \leq
b}|F_i|$. Take the union of the arrays arising from each row of
$C$ to obtain an FPA of size $N \mathrm{min}_{1 \leq i \leq
b}|F_i|$. Agreement between rows of this FPA corresponding to
different rows of $C$ can occur only at positions where the rows
of $C$ agree, since the symbol sets are disjoint. There are at
most $c$ such positions, so any two rows of the new FPA have
distance at least $c \geq bd$, and the array is an
$FPA_{\lambda}(bn,bd)$. $\hfill \square$

A useful tool in building new arrays from old is the direct
product.
\begin{Prop}\label{DirectProd}
Let $X_1$ be an $FPA_{\lambda}(n_1,a)$ of size $N_1$ and let $X_2$
be an $FPA_{\lambda}(n_2,b)$ of size $N_2$.  Then an
$FPA_{\lambda}(n_1+n_2,\mathrm{min}(a,b))$ may be constructed of
size $N_1 N_2$.

In particular, for even $n$, given two
$FPA_{\lambda}(n,\frac{n}{2})$, of sizes $N_1$ and $N_2$
respectively, an $FPA_{\lambda}(2n,\frac{n}{2})$ may be
constructed of size $N_1 N_2$, so that
\[ M_{\lambda}(2n,\frac{n}{2}) \geq
M_{\lambda}(n,\frac{n}{2})^2.\]

\end{Prop}
\noindent {\bf Proof}  Relabelling if necessary, construct $X_1$
and $X_2$ on disjoint symbol sets, giving
$\frac{n_1+n_2}{\lambda}$ symbols in total. Take the direct
product of $X_1$ and $X_2$, i.e. $Y=\{(u,v): u \in X_1, v \in X_2
\}$, where an ordered pair of codewords is interpreted as their
concatenation. Now, $Y$ is a set of $\lambda$-permutations of
length $n_1+n_2$, with frequency $\lambda$.  Any pair of
$\lambda$-permutations in $Y$ differ in at least
$\mathrm{min}(a,b)$ positions, hence $Y$ is an
$FPA_{\lambda}(n_1+n_2,\mathrm{min}(a,b))$.$\hfill \square$

In \cite{ColKloLin}, a permutation array is defined to be
\emph{$r$-separable} if it is a disjoint union of $r$ $PA(n,n)$'s
of size $n$.  We constructed an example of such a PA in part (i)
of Theorem \ref{UpperBound}.  We use this notion of a separable
PA, i.e. a PA which is a disjoint union of other PA's, in the next
result.

\begin{Thm}\label{SepPA}
\begin{itemize}
\item[(i)] Given a separable $PA(n,d)$ which is the disjoint union
of $r$ $PA(n,\delta)$'s, each of size $N$, where $2d \geq \delta$,
an $FPA_2(2n,\delta)$ of size $rN^2$ may be constructed.
\item[(ii)] Given $r$ MOLS of order $n$, an $FPA_2(2n,n)$ of size
$rn^2$ may be constructed.  If $n$ is a prime power, an
$FPA_2(2n,n)$ of size $(n-1)n^2$ is obtained.
\end{itemize}
\end{Thm}
\noindent {\bf Proof}  Denote the $r$ $PA(n,\delta)$'s by
$\Gamma_1,\ldots,\Gamma_r$.  For each $i=1,\ldots,r$, form the
direct product of $\Gamma_i$ with itself, i.e. $Z_i=\{(u,v): u,v
\in \Gamma_i \}$. Then $Z_i$ is a set of $N^2$
$\lambda$-permutations of length $2n$, on $n$ symbols, with
frequency $\lambda=2$, and minimum distance $\delta$.  Take the
union $Z=Z_1 \cup \ldots \cup Z_r$.  The $\lambda$-permutations
from different $Z_i$ have pairwise distance $2d \geq \delta$, and
hence $Z$ is an $FPA_2(2n,\delta)$ of size $rN^2$.

By a result established in \cite{ColKloLin} and reproved
constructively in \cite{HucMul}, $r$ MOLS of order $n$ may be used
to construct an $r$-separable $PA(n,n-1)$.  When used in the above
construction, this yields an FPA with $\delta=n$ and $2d=2n-2$
($>n$ for $n>2$), i.e. an $FPA_2(2n,n)$ of size $rn^2$. The last
part follows by noting that, for a prime power $n$, a complete set
of $n-1$ MOLS of order $n$ is obtainable. $\hfill \square$

\begin{Ex}\label{FPA284}
By the construction from part (ii) of Theorem \ref{SepPA}, an
$FPA_2(8,4)$ of size $48$ may be obtained from $3$ MOLS of order
$4$.  For example, the MOLS
\[ L_1=
\begin{array}{cccc}
0&1&2&3\\
1&0&3&2\\
2&3&0&1\\
3&2&1&0\\
\end{array},
L_2=
\begin{array}{cccc}
0&1&2&3\\
2&3&0&1\\
3&2&1&0\\
1&0&3&2\\
\end{array},
L_3=
\begin{array}{cccc}
0&1&2&3\\
3&2&1&0\\
1&0&3&2\\
2&3&0&1\\
\end{array}
\]

yield an FPA whose first $8$ rows are listed below.
$$
\begin{array}{cccccccc}
0&1&2&3&0&1&2&3\\
0&1&2&3&1&0&3&2\\
0&1&2&3&2&3&0&1\\
0&1&2&3&3&2&1&0\\
1&0&3&2&0&1&2&3\\
1&0&3&2&1&0&3&2\\
1&0&3&2&2&3&0&1\\
1&0&3&2&3&2&1&0
\end{array}
$$
\end{Ex}

The next two results generalize the direct product construction of
Theorem \ref{SepPA}.  A similar approach is explored in
\cite{ChuColDuk2}, in the context of constant composition codes;
the reader is referred to \cite{ChuColDuk2} for more details, and
for proofs of Theorem \ref{SepFPA} and Theorem \ref{SepFPA2}.

\begin{Thm}\label{SepFPA}
For $i=1,\ldots,b$, let $X_i$ be a separable
$FPA_{\lambda}(n,d_i)$ which is a disjoint union of $r_i$
$FPA_{\lambda}(n,\delta_i)$'s,
${\Gamma^i}_1,\ldots,{\Gamma^{i}}_{r_i}$, with $\sum d_i \geq
\mathrm{min}\{\delta_i\}$. Denote $\mathrm{min}\{\delta_i\}$ by
$\delta$, and $\mathrm{min}\{r_i\}$ by $r$.  Then an $FPA_{b
\lambda}(bn,\delta)$ may be constructed, of size $\sum_{j=1}^r
(\prod_{i=1}^{b} |{\Gamma_j}^{i}|)$.
\end{Thm}

The direct product construction in Theorem \ref{SepPA} and Theorem
\ref{SepFPA} may be adapted in various ways by choosing some
subset of the direct product which has special properties. In
\cite{ChuColDuk}, a recursive construction of PA's is given, which
uses transversal packings; the next result indicates one way in
which transversal packings may be used to construct an FPA from
separable PAs.  This construction may be applied to a set of
disjoint separable PA's such as those obtainable from the MOLS
construction of \cite{HucMul}, or to a single such PA with its
subarrays permuted appropriately.

\begin{Thm}\label{SepFPA2}
Let $X_1,\ldots,X_k$ be $k$ separable PA's, such that each
$X_i=PA(n,d_i)$ is a disjoint union of $r$
$PA(n,{d_i}^{\prime})$'s,
${\Gamma_i}^{(1)},\ldots,{\Gamma_i}^{(r)}$, and the
${\Gamma_i}^{(j)}$'s may be ordered such that
${\Gamma_1}^{(j)},\ldots,{\Gamma_k}^{(j)}$ are disjoint for each
$j$. Suppose there exist transversal packings $T_1,\ldots,T_r$ ,
where each $T_j$ has distance $\delta$ and type
$|{\Gamma_1}^{(j)}| \ldots |{\Gamma_k}^{(j)}|$.   Denote $d_1+
\cdots + d_k$ by $D$, and denote the smallest sum of any $\delta$
of the ${d_i}^{\prime}$ by $t$. Then an $FPA_k(kn,d)$ may be
constructed, of size $\sum_{j=1}^{r}|T_j|$, where $d=
\mathrm{min}(t,D)$.
\end{Thm}

We note that Theorem 3.2 of \cite{ChuColDuk} may also be
generalized to construct an $FPA_{\lambda}(n,d)$ from $k$
separable $FPA_{\lambda}(n_i,d_i)$'s ($1 \leq i \leq k$).
Replacing the $PA(n_i,d_i)$'s by the equivalent FPA's in the
statement and proof of this result, an immediate generalization
for $\lambda>1$ is obtained.

\section{Special cases}

An $FPA_{\lambda}(n,d)$, where $n=m \lambda$, may be viewed as an
$m$-ary code with constant weight composition $(\lambda, \ldots,
\lambda)$. In certain special cases, known results for constant
weight codes provide bounds and constructions of relevance to
FPAs.

\begin{Prop}
If $n=2 \lambda$, then an $FPA_{\lambda}(n,d)$ of size $M$ is a binary code $(n,M,d)$
of length $n$, minimum (Hammimg) distance $d$ and constant weight $\lambda$.
\end{Prop}

In \cite{Brouwer}, constructions and bounds are given for
$A(n,d,w)$, the maximum possible number of binary vectors of
length $n$, Hamming distance at least $d$ and constant weight $w$,
for values of $n$ up to $28$.  Observe that
$A(n,d,\frac{n}{2})=M_{\frac{n}{2}}(n,d)$.  The exact value of
$A(n,d,w)$, and corresponding constructions, is known for all
lengths $n \leq 11$.  If $d$ is odd, then
$M_{\frac{n}{2}}(n,d)=M_{\frac{n}{2}}(n,d+1)$, so only even
distances need be considered. The following FPAs may be directly
constructed, by use of Hadamard matrices and Steiner systems
(\cite{Brouwer}).

Recall that a \emph{Hadamard matrix} is a square matrix with
entries $+1$, $-1$ whose rows are mutually orthogonal.  Hadamard
matrices of order $n$ can only exist for $n=1,2$ and $n=4k$; it is
conjectured that they exist for each $n=4k$.  Hadamard matrix
constructions and properties may be found in Section IV.24 of
\cite{ColDin}.

\begin{Thm}[Theorem 10, \cite{Brouwer}] \label{HadThm}
$M_{\frac{n}{2}}(n,\frac{n}{2})=2n-2$ if and only if a Hadamard
matrix $H_n$ of order $n \geq 1$ exists.
\end{Thm}

An (optimal) $FPA_{\frac{n}{2}}(n,\frac{n}{2})$ may be constructed
from the Hadamard matrix $H_n$ as follows.  First convert the
entries of the `half-frame' of $+1$'s bordering $H_n$ into $-1$'s.
Now take the non-initial rows of $H_n$ and $-H_n$, and convert the
entries $+1$ to $0$ and $-1$ to $1$ in every row.

\begin{Ex}\label{HadMat} Using the Hadamard matrix of order $12$ (unique up to
isomorphism)
gives an $FPA_6(12,6)$ of size $22$.  We list the first few rows.

$$\begin{array}{cccccccccccc}
1&0&1&0&1&1&1&0&0&0&1&0\\
1&0&0&1&0&1&1&1&0&0&0&1\\
1&1&0&0&1&0&1&1&1&0&0&0\\
1&0&1&0&0&1&0&1&1&1&0&0
\end{array}
$$
\end{Ex}

Combining Theorem \ref{HadThm} with Proposition \ref{DirectProd}
we see that, if $n$ is an even number such that a Hadamard matrix
of order $n$ exists, then $M_{\frac{n}{2}}(2n,\frac{n}{2}) \geq
(2n-2)^2$.  For example, $M_6(24,6)>484$.

A \emph{Steiner system} $S(t,k,v)$ is a $t-(v,k,1)$ design, that
is, a collection of $k$-subsets (called blocks) of a $v$-set such
that each $t$-tuple of elements of this $v$-set is contained in a
unique block.   When $t=3$ and $k=4$, this called a \emph{Steiner
quadruple system}.

\begin{Ex}
Using the Steiner quadruple system $S(3,4,8)$, an $FPA_4(8,4)$ of
size $14$ may be constructed.   The extended cyclic code
$\{(1011000)1, (0100111)0 \}$ is one example; the code is
constructed by taking cyclic developments of the vectors in
parenthesis.
\end{Ex}

We conclude by remarking that, in the study of $PA(n,d)$ arrays,
one builds the rows of the array by using permutations on $n$
symbols, and in $FPA_\lambda(n,d)$ arrays, one builds rows by
using $m$ distinct symbols, each repeated exactly $\lambda$ times.
However, there is in fact no need for such uniformity of
frequency, and one could consider the following, very general,
setting.

Let $n=\lambda_1 +\cdots +\lambda_r$ be a partition of $n$. Then
one could consider constructing arrays with the property that in
each row, for $i=1, \dots, r$, the symbol $i$ occurs exactly
$\lambda_i$ times. From papers such as \cite{ChuColDuk}, there is
motivation for studying such a general setting; in fact the
corresponding constant composition codes have been widely studied;
see for example \cite{Brouwer}. Sets of $F(n; \lambda_1, \dots,
\lambda_r)$ orthogonal frequency squares have been studied (see
Chapter 4 in \cite{LayMul}).  However, we do not consider
frequency permutation arrays with an arbitrary frequency vector
$n= \lambda_1 + \cdots + \lambda_r$ in this paper.

{\bf Acknowledgment:} The first author is supported by a Royal
Society Dorothy Hodgkin Fellowship.  The authors wish to thank the
anonymous referees for their comments and suggestions.

School of Mathematics and Statistics, Mathematical Institute, North Haugh, St. Andrews, Fife, KY16 9SS, United Kingdom; Email: sophieh@mcs.st-andrews.ac.uk

Department of Mathematics, The Pennsylvania State University, University Park, PA 16802, U. S. A.; Email: mullen@math.psu.edu

\end{document}